\documentclass[12pt]{article}
\usepackage{amsmath}
\usepackage{amsfonts}
\usepackage{graphicx}
\usepackage{color}
\usepackage{graphics}
\usepackage{graphicx}
\usepackage{amsmath}
\usepackage{amsfonts}
\usepackage{amssymb}

\numberwithin{equation}{section}

\begin{document}

\noindent{\bf \Large Inverse medium problems, saddle point formulation}

\vspace{2mm}

\noindent{Kazufumi Ito, North Carolina State University, Raleigh, USA}

\vspace{2mm}

\noindent\underline{Abstract:}   In this paper we discuss  
inverse medium problems. We develop the direct sampling 
method based on probing indices using the saddle point 
formulation. The medium is constructed by solutions of saddle point problems. The method improves the probing functions
for the direct sampling method and directly images the medium.
The method is very efficient and can be applied to a general class of  inverse medium problems.

\section{Introduction}
In this paper we consider a general class of the inverse medium 
problems \cite{IJ} and develop the saddle point formulation to construct the medium.  In order to motivate our proposed method we discuss
the Schrödinger or DOT problem.  Diffusive optical tomography (DOT)
\cite{CILZ}  is a popular noninvasive imaging technique that measures the optical properties of a medium and creates images which show the distribution of absorption coefficient inside the body. Medical applications of DOT include breast cancer imaging. We consider the model:
\begin{equation} \label{DOT}
-\Delta u+\mu(x)\,u=0,\;\ \frac{\partial u}{\partial \nu}=g,\;\;u=f\;\;
\mbox{   at  the boundary $\Gamma=\partial\Omega$}.
\end{equation}
Here, $\mu=\mu(x)$ is unknown potential.  The problem is to determine
the medium $(u,\mu)$ form the Cauchy data $(f,g)$
at the boundary $\Gamma$ of domain $\Omega$. 
It is very ill-posed problem \cite{IJ} and they are not uniquely identifiable.
We define the induced source $\lambda=\mu u$. It provides us
the decomposition of difficulties.  We solve the inverse source 
problem 
for $\lambda$ by the probing index and then determine
$\mu(x)$. We develop the probing index for the inverse medium problem  in 
\cite{IJ,IJZ,CIZ1,CIZ2,CILZ} and Appendix.
The probing index $\lambda$ approximates
the induced source $\mu(x)\,u$ and satisfies
$$\begin{array}{l}
-\Delta u^0+\lambda^0=0, \frac{\partial u^0}{\partial \nu}=g,
\\ \\
-\Delta \lambda=0, \frac{\partial u^0}{\partial \nu}=\frac{P(f-u_0)}{\epsilon},
\end{array} $$
where $\lambda^0$ is the initialization.

\subsection{Iterative probing}
We improve the probing index $\lambda$ by
the iterative refinement:
$$\begin{array}{l}
-\Delta u^k+\lambda^k=0, \frac{\partial u^k}{\partial \nu}=g,
\\ \\
-\Delta \lambda^{k=1}=0, \frac{\partial \lambda^{k+1} }{\partial \nu}=\frac{P(f-u^k)}{\epsilon}
\end{array} $$
where it is the form of ADI (alternate direction iteration) and
$\lambda^k$ is a sequence of estimates of  induced source 
$\mu u$ and $u^k$ is a sequence of estimated states.

\subsection{Saddle point formulation}
The iterative probing converges to the saddle point problem for 
$(u,\lambda)$:
$$\begin{array}{l}
-\Delta u+\lambda=0, \frac{\partial u}{\partial \nu}=g,
\\ \\
-\Delta \lambda=0, \frac{\partial \lambda }{\partial \nu}=\frac{P(f-u)}{\epsilon}
\end{array} $$
where $\epsilon$ is a noise level for the data fitting $u=f$
and $P=-(\Delta)_\Gamma$ is the surface Laplacian.
Moreover as $\epsilon\to 0^+$ $u=u_\epsilon$ converges to 
$\bar{u}$ that satisfies the bi-harmonic equation
$$
\Delta^2 \bar{u}=0,\;\; \frac{\partial \lambda }{\partial \nu}=g,\;\;
u=f\mbox{  at }\Gamma.
$$

\noindent{\bf Remark (Bi-harmonic)}

(1) The saddle point formulation is equivalent to
the harmonic source assumption $\Delta \lambda=0$
on $\lambda$.
Then, $(\bar{u},\lambda)$ satisfies the saddle point problem
and the Cauchy data.
\vspace{2mm}

(2) We construct the medium by  $\mu(x)=\lambda(x)/u(x)$ from the induced source  assumption.

\section{General case}
In this section we discuss a general  inverse medium of the form
$$
Lu+F(u,\theta)=0, B_1u=g,\;\; B_2u=f
$$
where $L$ is a closed linear operator on a Hilbert space $X$ and $F(u,\theta)$
is a (non) linear parameter ($\theta$)-dependent operator.
WE assume $(Lu,u)>\delta\,|u|^2$ for some $\delta>0$.
$\lambda=F(u,\theta)$ defines the induced source in  $X$
The iterate probing is given by
\begin{equation} \label{prob}
\begin{array}{l}
Lu^k+\lambda^k=0, \;\;B_1u^k=g
\\ \\
L\lambda^{k+1}=0, \;\;B_2^*\lambda^{k+1}=\frac{P(f-u^k)}{\epsilon}
\end{array} \end{equation}
where we assume
$$
(u,L^*\lambda)=(B_2^*\lambda,u)-(B_1u)+(Lu,\lambda).
$$

Now, we show that the iterative probing method
\eqref{prob} is convergent. \vspace{2mm}

\noindent{Theorem 1} $(u_k,\lambda_k)$ converges to the saddle point solution satisfying
$$\begin{array}{l}
Lu+\lambda=0, \;\;B_1u=g
\\ \\
L^*\lambda=0, \;\;B_2^*\lambda=\frac{P(f-u)}{\epsilon}.
\end{array} $$

\noindent{Proof:} Since
$$\begin{array}{l}
L(u^{k+1}-u^k)+\lambda^{k+1}-\lambda_k=0, \;\;
B_1(u^{k+1}-u^k)=0
\\ \\
L^*(\lambda^{k+1}-\lambda^k),\;\; 
B_2^*(\lambda^{k+1}-\lambda^k)
=\frac{1}{\epsilon}P(u^{k}-u^{k-1}).
\end{array} $$
Multiplying $\lambda^{k+1}-\lambda^k$ to the first equation and 
$u^{k+1}-u^k$ to the second equation we obtain
$$
|\lambda^{k+1}-\lambda^k|^2
+\frac{1}{\epsilon}(u^{k+1}-u^{k},P(u^{k}-u^{k-1})=0
$$
Multiplying $u^{k+1}-u^k$ to the first equation
$$
(L(u^{k+1}-u^k),u^{k+1}-u^k)+(u^{k+1}-u^k,\lambda^{k+1}-\lambda^k)
$$
and multiplying $u^{k}-u^{k-1}$ to the second equation
$$
\frac{1}{\epsilon}(P(u^k-u^{k-1},u_k-u^{k-1})
+(\lambda^k,\lambda^{k-1})=0.
$$
Thus, there exits a $\rho<1$ such that
$$
|L(|u^{k+1}-u^k)|^2+|\lambda^{k+1}-\lambda^k|^2\le \rho\,
|L(|u^{k}-u^{k-1})|^2+|\lambda^{k}-\lambda^{k-1}|^2. \square
$$

The next theorem proves that $(u_\epsilon,\lambda_\epsilon)$
converges to the inverse medium pair $\bar{u},\bar{\lambda})$
satisfying
\begin{equation} \label{lmd}
L\bar{u}+\bar{\lambda}=0,\;\; B_1\bar{u}=g,\;\; B_2\bar{u}=f.
\end{equation}

\noindent{Theorem 2} As $\epsilon \to0$
$(u_\epsilon,\lambda_\epsilon)$ converges to $(\bar{u},\bar{\lambda})$ where $\bar{u}$ satisfies
$$
L^*L\bar{u}=0, \;\; B_1{\bar u}=g,\;\;B_2{\bar u}=f.
$$
Proof: Multiplying $\lambda$ to the first equation and 
$u-\bar{u}$ to the second equation we obtain
$$
|\lambda-\bar{\lambda}|^2+|\lambda|^2+
\frac{1}{\epsilon}(P(u-\bar{u}),u-\bar{u})=|\bar{\lambda}|^2
$$
Thus, $\lambda$ is bounded and $|\lambda_\epsilon
-\bar{\lambda}$ to $0$ and $u \to \bar{u}$.  $\square$

Moreover, $\bar{\lambda}$ is a minimum norm solution of \eqref{lmd}.   

\subsection{Least square and Constrained saddle point formulation}
Since in general $\lambda=F(u,\theta)$ defines constraints on
$\lambda$.
To this end we develop the constrained saddle point formulation.
Consider the least squares approach  \cite{CIY}
$$
\min\quad |Lu|^2+\frac{1}{\epsilon}(P(f-u),f-u)
$$
Equivalently, we have the constrained optimal control problem for source $\lambda$:
$$
\min \quad |\lambda|^2+\psi(\lambda)+\frac{1}{\epsilon}(P(f-u),f-u))
$$
subject to
$$
Lu+\lambda=0,\quad B_1u=g,
$$
where $\psi$ is the constraint cost for $\lambda$, e.g.,  the
sparsity and non-negativity.
The necessary optimality is given
$$\begin{array}{l}
Lu+p=0,\;\,p \in v+\partial \psi(\lambda),\;\; B_1u=g
\\ \\
L^*p=0, \;\; B_2p=\frac{P(f-u)}{\epsilon}
\end{array} $$
For example the non-negative constraint  $v\ge0$ we have 
$\lambda=\max(0,p)$.  It is the constrained saddle problem 
for $(u,p)$. In the case of \eqref{DOT} 
$$
\mu(x)=\max(0,p(x))/u(x)
$$
assuming $u(x)>0$.

\subsection{High order form}
We consider the high order form
$$
P \lambda =0 \mbox{  and } PL<0
$$
For example
$$
P=L^*+L^*L.
$$

If $Ly=-\Delta y$, then
$$
PL=\Delta \Delta (I-\Delta).
$$
We consider the hyperbolic case
$$
Lu=H\cdot\nabla u,\;\;L^*\lambda=\nabla \cdot(-H\lambda)
$$
Then, we have
$$
(L^*(LLu),u)=(L^*Lu,Lu)<0
$$

\subsection{Time-dependent Case}

Next, we consider the time-dependent case  of heat equation 
$\tilde{L}u=u_t+Lu$ and wave equation $\tilde{L}u=u_{tt}-Lu$.
$\tilde{L}^*$ is the backward heat equation and wave equation.
our approach  coincides with the time-reversal method.
Consider the backward heat equation, i.e. $u(T)$ is given and determine  the initial condition $u(0$.
$$\begin{array}{l}
u_t+Lu+\lambda=0
\\ \\
-\lambda_t+L^*\lambda=0\quad \lambda(T)=\frac{y(T)-z}{\epsilon}.
\end{array} $$
$$\begin{array}{l}
u_{tt}-Lu+\lambda=0
\\ \\
\lambda_{tt}+L^*\lambda=0\lambda(T)=\frac{y(T)-z}{\epsilon}.
\end{array} $$

Consider the least squares formulation \cite{CIY}
$$
\min |\frac{dy}{dt}u-Ly|^2+|\nabla y(0)|^2+\frac{|y(T)-z|^2}{\epsilon}
$$
Equivalently, we have the optimal control problem of the form
$$
\min \quad |u|^2+|\nabla y(0)|^2+\frac{1}{\epsilon} |y(T)-z|^2
$$
subject to
$$
\frac{dy}{dt}y-Ly+u=0
$$
Then, the necessary optimality is given by
$$\begin{array}{l}
y_t+Ly+\lambda=0,\;\lambda(0)=\nabla y(0)
\\ \\
-\lambda_t+L^*\lambda=0,\;\; \lambda(T)=\frac{y(T)-z}{\epsilon}.
\end{array} $$

\subsection{Side-way  heat equation}
We consider the side way heat equation
with the Cauchy data 
$(f(t),g(t))$ at the boundary. The saddle point problem is
formulated as
$$\begin{array}{l}
y_t+Ly+\lambda=0,\quad \frac{\partial u}{\partial \nu}=g(t)
\\ \\
\lambda_t+L^*\lambda=0,\quad \frac{\partial \lambda}{\partial \nu}=\frac{P(f(t)-u(t)}{\epsilon},
\end{array}$$
which is the system of the forward equations for $(y,\lambda)$.
Specifically, we consider the moving potential
$$
\lambda(t,x)\mbox{  corresponds to  } \chi_t(x-\gamma(t))u
$$
where $\chi_t$ is the characteristic function for the shape of  the inclusion and  $\gamma(t)$ 
ie the motion of the moving potential.

\section{Appendix, Direct Sampling and Probing method}

In this section we introduce the so-called Direct Sampling method.
Consider the inverse medium problem:
given a pair $(f,g)$ determine $u$ and the medium distribution 
$\theta$ such that
\begin{equation} \label{inv}
{\cal L}(u,\theta)=g,\;\; Cu=f,
\end{equation}
where $u\to {\cal L}(u,\theta)$ is  an elliptic operator and $\theta$ represents the medium function. This model  ${\cal L}(u,\theta))$
is used in the concrete examples including
EIT, DOT, Helmholtz, Elastic, Maxwell. Heat equation, Wave equation) as well as nonlinear problems in \cite{IJ,IJZ,CIZ1,CIZ2,CILZ}.

If the $(f,g)$, say Cauchy data is given,
then is an over determined Cauchy problem 
for $u$  if $\theta$ is known.
Thus, our objective is  to determine $(u,\theta)$ from  data $(f,g)$ so that \eqref{inv} holds.
It may not have a unique solution $\theta$ and we normally formulate the inverse problem as
the nonlinear optimization problem of the form
$$
\min\;\;|f-Cu|^2_\Gamma,\;\; \mbox{  over $\theta\in {\cal C}$},
$$
where $C$ is the observation map and $f$ is the corresponding observation and ${\cal C}$ is a constraint set to determine the medium  $\theta$.  We solve 
$$
{\cal L}(u,\theta)=g,
$$
given $(\theta,g)$ for $u$ then  minimizes the least square cost for $\theta$, i.e., it is the nonlinear constrained regression problem  for $\theta$.  Usually it is solved by the iterative method such as the gradient and Newton method, e.g,
\cite{IK,IJ} and It can be be very costly and slow convergent due to the ill-posedness of the inverse problem. 

Our proposed direct sampling method differs
from the optimization approach and develop a direct method to probe the medium 
$\theta$.
We shall first discuss the general philosophy of DSMs. A DSM was introduced and studied in using near-field data for locating inhomogeneities in inverse acoustic medium scattering. 
The method provides the location of the inhomogeneities based on an indicator function, which is defined as the $L^2$ inner product over the measurement surface between the measured data and the fundamental solution of the Helmholtz equation. Numerical experiments have shown that this method is effective and very robust to noise. It is particularly successful in locating multiple clustered objects sitting inside an acoustic medium with only a single or a few incident plane waves. The success of the method comes from the following two observations: the scattered field can be approximated by a finite sum of the fundamental solution of the Helmholtz equation centered at inhomogeneous inclusions; and the fundamental solutions centered at different points are nearly orthogonal in the $L^2$ inner product over the measurement surface. These observations motivate our current development of DSMs to the other inverse medium problem, when only one pair of Cauchy data is available. We shall construct an appropriate family of probing functions such that the boundary data can be approximated as a linear combination of probing functions. Given a nominal parameter
$\bar{\theta}$ we select the probing factions based on the fundamental solution $\eta_x$ to
$$
{\cal L}(\eta_x,\bar{\theta})=\delta_x
$$
We  formulate the direct sampling  method to probe the medium by the testing the observation data 
against the probing functions. It results in the probing index that depict the distribution of the medium.
It is a form of migrating  the boundary data to interior points $x$ by the back-projection.

We introduce the direct sampling method  in terms of the inverse medium problem for the diffusive optical tomography (DOT):
$$
-\Delta u+\mu(x)u=0,\;\; \frac{\partial u}{\partial \nu}=g(x),\;\ u=f \mbox{  at the boundary $\partial\Omega)$}.
$$
Problem is  to determine the absorption coefficient $\mu(x)\ge0$ and $u$ from
a single Cauchy data $(f,g)$. 

We introduce the basic concept, First, we have the solution representation
$$
u(\xi)-u_0(\xi)=\int_\Omega G(\xi,y)f(y)\,dy=Af
$$
where $G$ is the Green's function satisfying
\begin{equation} \label{Type I}
-\Delta G+\bar{\mu}G=\delta_x,\;\;\frac{\partial G}{\partial \nu}_\Gamma=0.
\end{equation}
Next, if $f$ is the induced source function defined by
$$
f(y)=(\mu-\bar{\mu})(y)u(y)
$$
By the quadrature rule: for $\xi \in\Gamma$
\begin{equation} \label{Green}
u(\xi)-u_0(\xi) \sim \sum_i G(\xi,y_i)f(y_i)vol(V_i)=A^Nf^N,
\end{equation}
where $y_i$ is the center of volume $V_i$. That is, in this case we try to determine
the sum $f^N$ of points sources from  the boundary measurement $u-u_0$ at $\Gamma=\partial\Omega$.  Assume that the probing function $\eta_x(\xi)$ satisfies
$$
u(\xi)-u_0(\xi) \sim \sum_i\; \eta_{x_i}f(x_i)
$$
and
$$
(\eta_{x_i},\eta_{x_j})_{\Gamma}\mbox{  is diagonally dominant }
$$
I.e.
$$
\eta_{x_i} \mbox{   are nearly orthogonal.}
$$
Then, we have
$$
(u(\xi)-u_0(\xi,\eta_i))_\Gamma \mbox{  approximates  a distribution $f_i=f(x_i)$}
$$
It follows from this observation that
we define the probing index $I(x)$ by the  wighted Integral of  measurement data $u-u_0$ at the boundary against  the probing function
$$
\eta_x=G(\xi,x)|_\Gamma 
$$
at the boundary $\Gamma$, I.e.
$$
I(x)=\frac{(\Lambda \eta_x,u-u_0)_\Gamma}{|\eta_x|_W},\;\; x\in \Omega,
$$
where  $\Lambda$ is a (differential) operator on $\Gamma$ and $|\eta_x|_W$ is a scaling factor. A proper choice of Sobolev scales for $\Lambda$ and $W$  to enhances the probing.
That is, $W$  removes the boundary layer of the index.
The weighted inner product defines the similarity of $\eta_x$ and $u-u_0$ so that $I(x)$ defines the estimate of the distribution of $f$. 
Here the probing $I(x)$ represents an approximated distribution of  $f$, that is
$$
I(x)=A^*P(u-u_0)|_\Gamma \sim f(x).
$$ 
where $P$  id the scaling defined by $\Lambda$ and $W$.
A proper choice  Sobolev scales for $\Lambda$ and $W$   enhances the probing index substantially.
 
Define  the integral operator  (the kernel function) $C$ from $L^2(\Omega)$
to $L^2(\Omega)$  by
$$
C(x,y)=\frac{(P \eta_x,\eta _y)_\Gamma}{|\eta_x|_W}.
$$
If the integral operator  is invertible, then we have  the identifiability of $f$ from the boundary measurement of $u-u_0$. If $C$ is diagonally dominant then the probing function $\eta_x$ is nearly orthogonal. The discretized version $C$ is given by
$$
C^N=(A^N)^*P^NA^N
$$

\noindent\underline{\bf Remark}

\noindent (1) Since the probing function is analytic at the boundary  it is less sensitive to  noise in data $y$,
i.e., the averaging theorem is applied.
\vspace{2mm}

\noindent (2) One can project $I(x)$ to the constraint  set such as the nonnegative constraint.
\vspace{2mm}
 
\noindent (3) The least square solution (maximum likelihood)  
to $Ax-f=0$ is given by
$$
x=A^*(AA^*)^{-1}f
$$
and it is significantly ill-posed.  
\vspace{2mm}

\noindent (4) If $P=I$ then $A^*f$ defines the back projection of $f$ and $A^*Pf$ defines the filtered back projection where 
$P$ represents the filtered version of $(AA^*)^{-1}$.  It is a direct method to the determine the  induced souse distribution
$$
I=A^*Pf.
$$

\noindent(5) $I(x)=\lambda=A^*Pf$ solves the adjoint equation,
$$
-\Delta \lambda+\bar{\mu}\lambda=0,\;\;\frac{\partial \lambda}{\partial \nu}=P(f-u).
$$
since
$$
(-\Delta \lambda+\bar{\mu}\lambda,\eta_x)
=-(\frac{\partial\lambda}{\partial n},\eta_x)+(\frac{\partial\eta_x}{\partial n},\lambda)
+(-\Delta \eta_x+\bar{\mu}\eta_x)=\lambda(x)
$$
So, we do not need to construct the fundamental solution $\eta_x$
to evaluate the index $I$.
\vspace{2mm}

\noindent(6) One can use the Drichlet-to-Newman formulation, i.e.
Assume $u=f$  as an input and use the observation 
$\frac{\partial u}{\partial n}=g$  to determine $\mu$ by
$$\begin{array}{l}
-\Delta u+\bar{\mu}u+\lambda=0,\;\;  u=f \mbox{  at  } \Gamma
\\ \\
-\Delta \lambda+\bar{\mu}\lambda=0,\;\;  
\frac{\partial}{\partial n}\lambda=\frac{1}{\epsilon}
P(g-\frac{\partial}{\partial n}u_0) \mbox{  at  } \Gamma.
\end{array}$$

\end{document}